\documentclass{icmart}

\contact[p.m.topping@warwick.ac.uk]{Mathematics Institute, University of Warwick, Coventry,
CV4 7AL, UK}

\usepackage{mathrsfs}
\usepackage{amsfonts, amsmath, wasysym}
\usepackage{amssymb,paralist} %amsthm didn't work here.
\usepackage{latexsym,nicefrac}
\usepackage[usenames]{color}
\usepackage{url}

\ifx\pdfoutput\undefined
\usepackage{graphicx} 
\else 
\usepackage[pdftex]{graphicx} 
\usepackage{epstopdf} 
\fi

\usepackage{tikz}

\newtheorem{theorem}{Theorem}[section]

\theoremstyle{definition}

\newtheorem{remark}[theorem]{Remark}
\newtheorem{example}[theorem]{Example}

\newcommand{\mm}{\ensuremath{{\cal M}}}
\newcommand{\pl}[2]{{\frac{\partial #1}{\partial #2}}}
\newcommand{\ga}{\gamma}

\newcommand{\vph}{\varphi}

\newcommand{\ep}{\varepsilon}

\newcommand{\RR}{\ensuremath{{\mathbb R}}}
\newcommand{\downto}{\downarrow}
\newcommand{\upto}{\uparrow}

\newcommand{\lap}{\Delta}

\DeclareMathOperator{\Vol}{vol}

\newcommand{\beq}{\begin{equation}}
\newcommand{\eeq}{\end{equation}}
\newcommand{\beqs}{\begin{equation}}
\newcommand{\eeqs}{\end{equation}}
\newcommand{\beqa}{\begin{equation}\begin{aligned}}
\newcommand{\eeqa}{\end{aligned}\end{equation}}
\newcommand{\beqas}{\begin{equation}\begin{aligned}}
\newcommand{\eeqas}{\end{aligned}\end{equation}}
\newcommand{\brmk}{\begin{remark}}
\newcommand{\ermk}{\end{remark}}
\newcommand{\partref}[1]{\hbox{(\csname @roman\endcsname{\ref{#1}})}}

\newcommand{\Rm}{{\mathrm{Rm}}}
\newcommand{\Ric}{{\mathrm{Ric}}}

\newcommand*\dz{\mathrm{d}z}

\title[Ricci flows with unbounded curvature]{Ricci flows with unbounded curvature\footnote{Proceedings of the International Congress of Mathematicians, Seoul 2014.}}

\author[Peter M. Topping]{Peter M. Topping}
\begin{document}

\begin{abstract}
Until recently, Ricci flow 
was viewed almost exclusively 
as a way of deforming Riemannian metrics of bounded curvature. Unfortunately, the bounded curvature hypothesis is unnatural for many applications, but is hard to drop because so many new phenomena can occur in the general case.
This article surveys some of the theory from the past few years that has sought to rectify the situation in different ways.
%
%The abstract should be comprehensible to any mathematician.
%A rough statement in words is preferable to a precise statement
%loaded with symbols and technical notions.
%There should be no references in the abstract.
%If a paper \textit{must} be mentioned, spell it out in full.
\end{abstract}



\maketitle

\section{Introduction}

Since its inception in 1982 \cite{ham3D}, Ricci flow has supported the development of a remarkable and elegant theory.
The flow has become well-known as a way of deforming a Riemannian metric in order to improve it, or turn it into a special metric that might satisfy a geometrically rigid condition or simply a natural PDE, and indeed up until now, most applications fit within the following general strategy.
\begin{itemize}
\item
First we take a space that we do not understand very well, perhaps a Riemannian manifold satisfying a curvature condition, or a metric space with some weak geometric structure. 
\item
Next we deform the space under Ricci flow, keeping track of its properties, for example its topology, its curvature or its conformal structure, until it develops into a very special space, for example sometimes one of constant curvature.
\item
Such special spaces we can hope to classify, and if the Ricci flow can be sufficiently well understood then we can go back and classify or better understand the space with which we started. 
\end{itemize}
As an example, Hamilton's original insight was that a simply connected three-dimensional closed Riemannian manifold  of positive Ricci curvature will flow smoothly under a suitably normalised Ricci flow through a family of manifolds of positive Ricci curvature to a manifold of \emph{constant} curvature, which can then be identified as a round sphere. He deduced that 
the original manifold of positive Ricci curvature must be diffeomorphic to the three-sphere. Dramatic further development of Ricci flow theory by Hamilton and then Perelman ultimately extended this principle to handle all closed three-manifolds, leading to a resolution of the conjectures of Poincar\'e and Thurston (\cite{formations, P1, P2, P3, KLnotes}).

This general strategy has clearly been very effective, but its scope has nevertheless been severely limited by the Ricci flow existence and asymptotics theory only applying in very special cases, meaning that we are only scratching the surface of the potential applicability of the method. In particular, the vast majority of applications require the underlying manifold to be closed, and without that hypothesis we are not even able in general to start the flow going, even for a short time, without imposing further conditions such as boundedness of the curvature that may damage potential applications.

This article surveys part of the programme to extend the theory of Ricci flow to handle general manifolds or even metric spaces. The central point will be the necessity to handle flows with unbounded curvature; until recently we have had no idea how to even start the flow going starting with a manifold of unbounded curvature, let alone understand its long-time existence and asymptotics or uniqueness. In this unrestricted situation, the flow interacts with itself `at spatial infinity' in an unfamiliar way that is interesting both geometrically and in terms of the pure PDE questions it raises.
Since classical solutions to the Ricci flow are characterised as existing until the curvature becomes unbounded, and we want to consider unbounded curvature from the outset, we now need to understand the issue of long-time existence in much more detail. This also brings to the fore the subtle issue of asymptotics of the flow, and we will witness how Ricci flow organises itself to find a special metric, even when there appear to be obstructions.

\section{Why is Ricci flow with unbounded curvature so difficult?}

We call a smooth one-parameter family of Riemannian metrics $g(t)$ on a manifold $\mm$ a Ricci flow when
\beq
\label{RFeq}
\pl{g}{t}=-2\,\Ric_{g(t)}
\eeq
where $\Ric$ is the Ricci tensor (see \cite{ham3D, RFnotes}).
One should interpret the right-hand side of the equation as some sort of Laplacian of the metric, and thus interpret this equation as some sort of heat equation, although ultimately this is a nonlinear equation, and is even not quite parabolic, which causes problems when trying to establish existence of solutions even to this day, except in relatively simple situations.

One way of trying to pose this flow is to consider an initial Riemannian metric $g_0$, and then look for a family $g(t)$, $t\in [0,T)$ satisfying \eqref{RFeq}, with $g(0)=g_0$. In due course, we will see that this is rather naive in general, but it works well in special situations such as on closed manifolds (Hamilton \cite{ham3D}) or more generally when the initial metric is complete and of bounded curvature (Shi \cite{shi}). The following hybrid of their work also incorporates a uniqueness assertion of Chen-Zhu \cite{chenzhu}; the proofs have been clarified and simplified by DeTurck/Hamilton \cite{deturck, formations} and Kotschwar \cite{kotschwar}.
\begin{theorem}[{\bf The Hamilton-Shi flow}]
\label{hamilton_shi}
Given a smooth, complete Riemannian manifold $(\mm,g_0)$ of \emph{bounded curvature}, 
there exist a unique $T\in (0,\infty]$ and Ricci flow $g(t)$ for $t\in [0,T)$ satisfying the equation
\eqref{RFeq}, the initial condition $g(0)=g_0$, and the properties that the curvature remains bounded for $t\in [0,T_0]$, for any $T_0\in [0,T)$, and that
if $T<\infty$ then 
$$\sup_\mm|\Rm_{g(t)}|\to\infty\quad\text{ as }t\upto T,$$
where $\Rm_{g(t)}$ is the curvature tensor of $g(t)$.
Moreover, $(\mm,g(t))$ is complete for all $t\in [0,T)$.
\end{theorem}
We will call the Ricci flow whose existence is asserted by
this theorem the \emph{Hamilton-Shi} Ricci flow.

\begin{figure}
\begin{center}
\def\svgwidth{300pt}
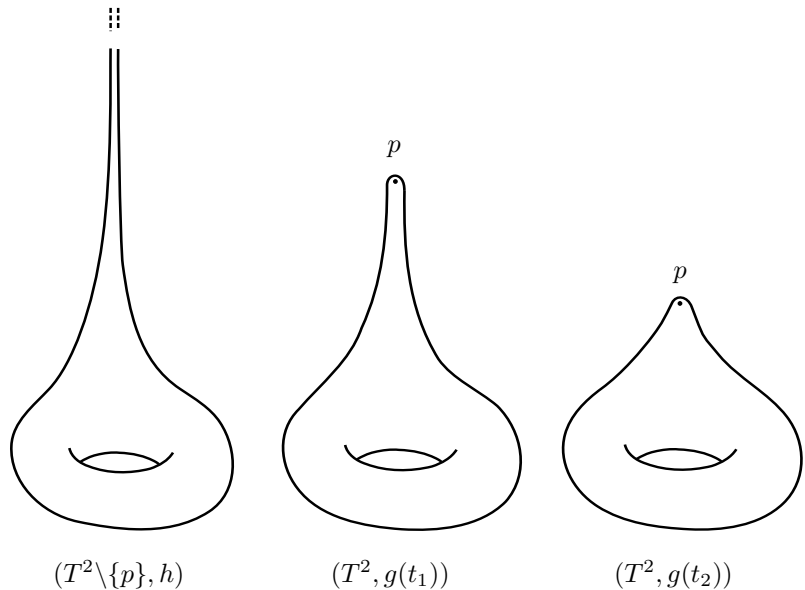
\end{center}
\caption{Pulling in points at infinity in Theorem 
\ref{revcusp_special}}
%%ALWAYS PUT LABEL AFTER CAPTION!!!!
\label{cusp_collapse_fig}
\end{figure}

The final assertion here that the Ricci flow is complete is more or less obvious. Indeed completeness of a Riemannian manifold is equivalent to the assertion that the length of any smooth proper curve $\ga:[0,1)\to\mm$ (i.e. any curve heading off to infinity in $\mm$) is infinite, but the boundedness of the curvature in this theorem guarantees that lengths of curves can only grow or decrease at most exponentially (see, for example, \cite[Lemma 5.3.2]{RFnotes})
and so an infinitely long curve at time $t=0$ will remain infinitely long at a later time $t\in [0,T_0]$. However, this principle emphatically fails for flows of unbounded curvature. Loosely speaking,
\begin{quote}
\emph{
Unbounded curvature allows Ricci flow to feel spatial infinity.
}
\end{quote}
In particular, an unbounded curvature Ricci flow can pull `points at infinity' to within a finite distance in finite time, as we now illustrate.
\begin{theorem}[{\bf Pulling in points at infinity;} Special case of \cite{revcusp}]
\label{revcusp_special}
There exists a smooth Ricci flow for $t\in [0,\infty)$, starting at a smooth, complete Riemannian metric of bounded curvature, that is incomplete for all $t>0$.

In particular, if $T^2$ is a torus equipped with an arbitrary conformal structure, $p\in T^2$ is any point, and we write $h$ for the unique complete, conformal, hyperbolic metric on $T^2\backslash \{p\}$, then there exists a smooth Ricci flow $g(t)$ on $T^2$ for $t> 0$ such that $g(t)\to h$ smoothly locally on $T^2\backslash \{p\}$ as $t\downto 0$.
\end{theorem}
Of course, the Ricci flow $g(t)$ on $T^2$ can be restricted to $T^2\backslash \{p\}$ to give the desired Ricci flow that starts complete, but instantaneously becomes incomplete. Intuitively the point $p$ is being pulled in from infinity as $t$ lifts off from zero -- see Figure \ref{cusp_collapse_fig}.

%\vskip 5mm
%\vskip 5mm

This example also illustrates the subtleties of uniqueness in Ricci flow in the presence of possibly unbounded curvature, since we could have flowed the initial manifold $(T^2\backslash \{p\},h)$ using the Hamilton-Shi flow instead. In this example that flow would simply  dilate the metric, and could be written explicitly as $g(t)=(1+2t)h$.

Examples such as this run somewhat counter to the classical intuition in Ricci flow, and the failure to recognise that they can arise in practice has previously led to errors in important parts of the literature -- see the discussion in \cite{HCT}.
We will return later to see how incompleteness can be dealt with, and also how to impose a condition at infinity (analogous to a boundary condition) in order to make the problem well posed.

\section{A clear picture of Ricci flow on surfaces}
\label{2Dsect}

Both the nature of the evolution equation and the existence theory above become crystal clear when the dimension of the underlying manifold $\mm$ is two.

In this case, the Ricci tensor can be written in terms of the Gauss curvature $K$ as $\Ric=K.g$, so the flow equation is
$$\pl{}{t}g(t)=-2K. g(t).$$
Thus %in this case, 
the flow will always make a conformal deformation of the metric, and we may take local isothermal coordinates
$x$ and $y$, and write the flow 
$g(t)=e^{2u}|\dz|^2:=e^{2u}(dx^2+dy^2)$
for some locally-defined, scalar, time-dependent function $u$, which can then be shown to satisfy the local equation
\beq
\label{2DRFeq}
\pl{u}{t}=e^{-2u}\lap u = -K.
\eeq
In particular, in this form, we are dealing with a strictly parabolic PDE. 
Up to a change of variables, this is the so-called logarithmic fast diffusion equation, which has an extensive literature -- see \cite{vaz, DK} and the references therein.

\begin{theorem}[{\bf Closed surfaces.} {Hamilton \cite{Ham88}, Chow \cite{chow}}]
\label{HCthm}
Let $\mm$ be a closed, oriented surface and $g_0$ any smooth metric. 
Define $T=\infty$ unless $\mm=S^2$ (topologically) in which case we set
$T=\frac{\Vol_{g_0}\mm}{8\pi}$.
Then there exists a unique Ricci flow $g(t)$ on $\mm$ for $t\in [0,T)$ 
so that $g(0)=g_0$. 
Depending on the genus of $\mm$, we have
%
%NOTE MANUAL SPACING
\begin{equation*}
\begin{array}{lll}
\mm=S^2:&
{\displaystyle \frac{g(t)}{2(T-t)}\to G_{+1},}&
\text{ a metric of const. curvature }+1,\text{ as }t\upto T.\\[12pt]
\mm=T^2:&
{\displaystyle g(t)\to G_{0},}&
\text{ a flat metric, as }t\to\infty.\\[5pt]
\mm\neq S^2, T^2:&
{\displaystyle \frac{g(t)}{2t}\to G_{-1},}&
\text{ a metric of const. curvature }-1,\text{ as }t\to\infty.
\end{array}
\end{equation*}
\end{theorem}
The well-posedness theory on closed surfaces was thus completed in the 1980s. An alternative approach to Theorem \ref{HCthm} that provides a model for the trickier argument required to prove the Poincar\'e conjecture can be found in \cite{park_city}.

It is apparent from Theorem \ref{HCthm} that Ricci flow \emph{geometrises} a closed surface; more precisely it finds a conformal metric of constant curvature on an arbitrary closed Riemann surface, which is enough to establish the Uniformisation theorem in the restricted case of closed surfaces \cite{CLT}.
It is then a natural question to ask whether Ricci flow performs the same geometrisation task on a general surface.
Of course, to be able to even ask this question, we have to be able to start the Ricci flow with a more general metric, and continue it until the flow has had a chance to organise itself into a special metric, whereas the Hamilton-Shi flow from Theorem \ref{hamilton_shi} flows restricted metrics, and in general will stop (as we will see) before the flow has achieved anything.

In fact, as 
we shall demonstrate shortly, 
Ricci flow is perfectly capable of flowing a completely general surface with possibly unbounded curvature in a uniquely defined way, without even requiring the initial metric to be complete.
Before stating the result, we dwell on some issues that such a result must address.

Those unfamiliar with PDE theory often misinterpret existence theory as presumably being obvious. Surely if we are trying to solve an equation $\pl{g}{t}=-2Kg$, then we should simply keep moving $g$ by tiny amounts in the direction $-2Kg$ and a solution should result. The most naive aspect of that suggestion is that it would appear to apply to the problem of solving backwards in time from a given smooth metric, whereas this is certainly not possible. Indeed, Ricci flow has the dramatic smoothing effect of parabolic equations, and immediately makes any metric real analytic \cite{kotschwar_ra}. Therefore we could never flow backwards in time starting from a general smooth metric that is not also real analytic.

A more subtle issue that arises once we drop the hypothesis that the underlying manifold is closed, is that we have to worry about boundary conditions. To illustrate this issue, consider simply the problem of starting the Ricci flow on the open disc $D^2\subset\RR^2$ with a metric $g_0=e^{2u_0}(dx^2+dy^2)$ that is smooth up to the boundary $\partial D^2$. From \eqref{2DRFeq}, this results in the PDE problem
\begin{equation}
\left\{
\begin{aligned}
\pl{u}{t} &=e^{-2u}\lap u &  \qquad & D^2\times [0,T)\\
u(0)&=u_0  & \qquad & D^2.
\end{aligned}
\right.
\end{equation}
Even amongst solutions that remain continuous up to the boundary at later times, this PDE problem is ill-posed owing to nonuniqueness. Indeed, standard parabolic theory tells us that we are free to specify the restriction of $u$ to $\partial D^2$ at later times, and only then would we obtain uniqueness. In fact, we will shortly solve the well-posedness problem without resorting to specifying any more data.

A third way we can see that existence theory should be a subtle issue is to consider what it should imply. Geometric flows are typically designed to find special objects, for example constant curvature metrics in the present discussion. On the other hand, it is often possible to give a geometric flow initial data that for some reason cannot be deformed globally and smoothly into a special object -- occasionally one does not even exist. Any assertion of global existence of solutions in such cases is also asserting that the geometric flow must organise itself in such a way to resolve this issue, often finding an ingenious way of decomposing the object being flowed into multiple special objects (see \cite{rick_survey} and \cite{RTZ, THMF_negcurv} for some other contexts in which this occurs).

We shall shortly see how these issues are resolved in practice by the Ricci flow, but first we give the main well-posedness result in the two-dimensional case. The existence part is joint work with G. Giesen.

\begin{theorem}[{\bf Ricci flows on surfaces, without restriction;} \cite{GT2, JEMS, ICRF_uniq}]
\label{fullthm}
Let $(\mm,g_0)$ be any smooth (connected) Riemannian surface, possibly incomplete and/or with unbounded curvature. 
Depending on the conformal type of $\mm$, we define $T\in(0,\infty]$ 
by\footnote{In the case that $\mm=\mathbb C$, we set $T=\infty$ if $\Vol_{g_0}\mathbb C=\infty$. Here $\chi(\mm)$ is the Euler characteristic.}
$$ T := \begin{cases}
    \frac1{4\pi\chi(\mm)}\Vol_{g_0}\mm & \text{if }(\mm,g_0)\cong\mathcal S^2, \mathbb C\text{ or }\mathbb R\!P^2, \\
    \qquad\infty & \text{otherwise}.
\end{cases} $$
Then there exists a 
smooth Ricci flow $g(t)$ on $\mm$, defined for $t\in[0,T)$ such that 
\begin{compactenum}
\item $g(0)=g_0$, and 
\item $g(t)$ is instantaneously complete, i.e. complete for all $t\neq 0$ at which it is defined.
\end{compactenum}
The flow $g(t)$ is unique in the sense that if $\tilde g(t)$ is any other smooth Ricci flow, defined now for $t\in [0,\tilde T)$, satisfying (1) and (2) above, then $\tilde T\leq T$ and $g(t)=\tilde g(t)$ for all $t\in [0,\tilde T)$.

In addition, the Ricci flow $g(t)$ is maximally stretched (see Remark \ref{max_stretch}), and the Gauss curvature $K_{g(t)}$ satisfies
$$K_{g(t)}\geq -\frac{1}{2t}$$
for $t\in (0,T)$.
If $T<\infty$ then we have 
\[
\Vol_{g(t)}\mm = 4\pi\chi(\mm) (T-t)\ \to\  0 \quad\text{ as } t\upto T, \] 
and in particular, $T$ is the maximal existence time.
\end{theorem}

Related results can be found in the literature of the logarithmic fast diffusion literature (e.g. \cite{DK, vaz}) and the work of Mazzeo, Sesum, Ji and Isenberg \cite{JMS09, IMS13}.

\begin{remark}[{\cite{GT2} and \cite[Remark 1.5]{GT4}}]
\label{max_stretch}
The \emph{maximally stretched}  assertion of the theorem means that $g(t)$ lies `above' any another Ricci flow with the same or lower initial data.
More precisely,
if $0\leq a<b\leq T$ and $\tilde g(t)$ is any Ricci flow on $\mm$ for $t\in [a,b)$ with $\tilde g(a)\leq g(a)$ (with $\tilde g(t)$ not necessarily complete or of bounded curvature) then $\tilde g(t)\leq g(t)$ for every $t\in[a,b)$. 
\end{remark}

The mechanism by which this Ricci flow makes an incomplete metric immediately complete is somewhat similar to how the example from Theorem \ref{revcusp_special} made a complete metric immediately incomplete. Unbounded curvature allows points at a finite distance to be sent out to infinity, and vice versa.

Given the unusual nature of these flows, some examples are in order. Although Theorem \ref{fullthm} can handle arbitrary initial metrics, with arbitrarily wild behaviour at spatial infinity, we pick the two simplest possible examples \cite{JEMS}, both of which start completely flat but have no choice but to have unbounded curvature immediately.

\begin{example}
\label{punc_plane_ex}
Let $(\mm, g_0)$ be the punctured plane. The metric is incomplete since we can take a curve asymptoting to the origin, of finite length, and the flow will have to do something about this in order to make the metric complete immediately.
The flow we construct in Theorem \ref{fullthm} deals with this by
stretching the metric near the puncture to be asymptotically a hyperbolic cusp, scaled to have Gauss curvature $-\frac{1}{2t}$.
See Figure \ref{punc_plane_fig}.
(In fact, in this extremely special case, one could compute the flow explicitly as a so-called Ricci soliton, up to the solution of an ODE.)

%\vskip 5mm
%%%%%NOTE: I HAVE RESTRICTED THE POSITIONING OF THE FOLLOWING FIGURE BY USING [h]
\begin{figure}[h]
\begin{center}
\def\svgwidth{250pt}
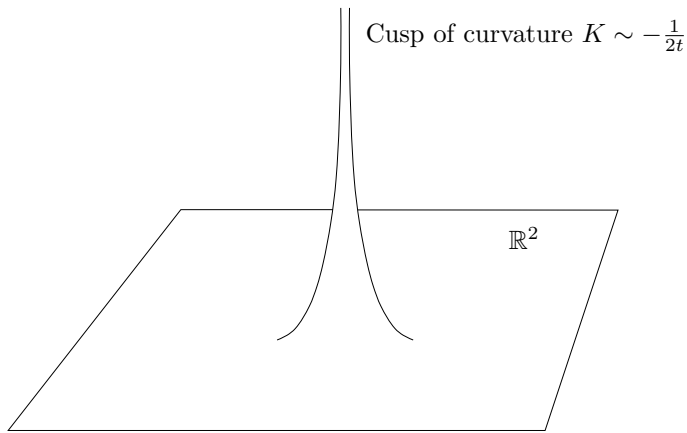
\end{center}
\caption{Puncture turns into a hyperbolic cusp in Example \ref{punc_plane_ex}}
\label{punc_plane_fig}
\end{figure}
%\vskip 5mm
\end{example}

\begin{example}
\label{disc_ex}
Now let $(\mm,g_0)$ be the Euclidean unit two-dimensional disc.
Again, the flow must blow up the metric to make it complete immediately, and it does this by stretching it near to the `boundary' circle into a Poincar\'e metric, also  scaled to have Gauss curvature $-\frac{1}{2t}$. See Figure \ref{disc_fig}.

%\vskip 5mm
\begin{figure}
\begin{center}
\def\svgwidth{250pt}
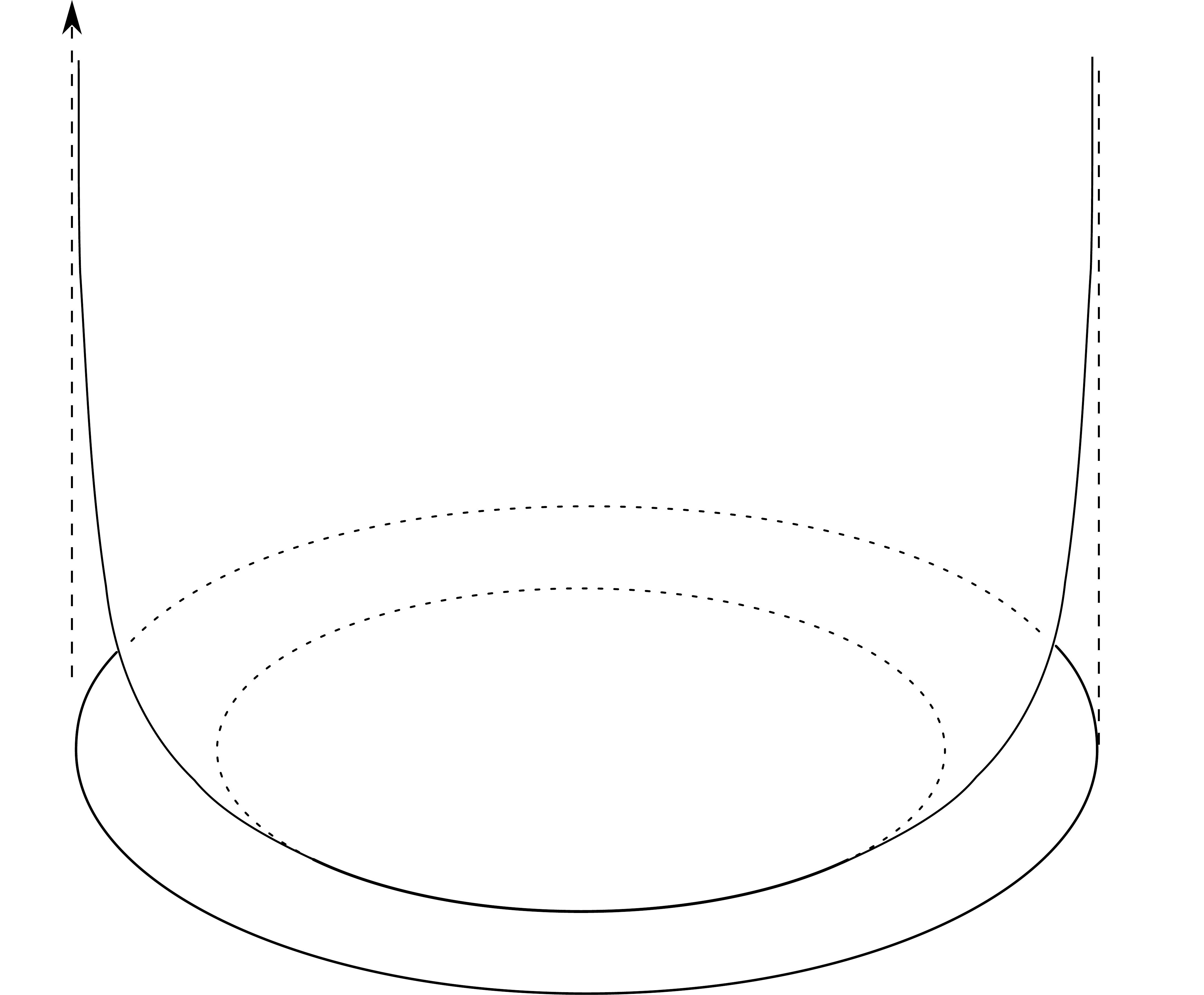
\end{center}
\caption{Metric stretches at infinity in Example \ref{disc_ex}}
\label{disc_fig}
\end{figure}
%\vskip 5mm

This example fits directly into the discussion above on flowing smooth metrics on the disc $D^2$, where we decided that the addition of boundary data would be the standard way of achieving well-posedness.
In Theorem \ref{fullthm}, the condition on solutions that they are instantaneously complete can be viewed as a substitute for a boundary condition.
\end{example}

Intuitively, the flow from Theorem \ref{fullthm} is finding a way of feeding in volume at infinity where the metric is incomplete. The notable feature here is not just that this can be done in order to give a solution, but that it can be done uniquely: If we try to feed in less volume, then the metric will fail to become complete.  On the other hand, if we try to feed in more, then the damping within the equation, arising from the $e^{-2u}$ factor in the equation \eqref{2DRFeq} is preventing the extra volume from arriving in the interior.

Naively, one might view a common feature of both these examples to be that the conformal factor $u$ in the most obvious, Euclidean coordinates is immediately asymptotically infinity at spatial infinity. However, the conformal factor $u$ depends on the coordinates chosen, and in different coordinates this property will not hold.

Given that Theorem \ref{fullthm} makes the metric complete immediately, and also makes the curvature bounded from below, it is reasonable to speculate that maybe the flow also makes the curvature bounded from above immediately, and that therefore after an arbitrarily short time we are in the classical situation of Theorem \ref{hamilton_shi} and could make do from then on with the Hamilton-Shi flow. 
We will see that this suggestion is wrong in two ways. 
To begin with, we cannot hope ever to be able to construct a complete Ricci flow solution that makes the curvature bounded from above, because our flow is the unique instantaneously complete solution, and carefully chosen initial metrics will have unbounded curvature for all time:

\begin{theorem}[{\bf Ricci flows with unbounded curvature;} \cite{GT3}]
\label{unbdd_ex}
Given any noncompact Riemann surface $\mm$, there exists a smooth, complete, conformal metric $g_0$ such that the unique complete Ricci flow $g(t)$ given by Theorem \ref{fullthm} exists for all $t\geq 0$ and $(\mm,g(t))$ has unbounded curvature for each $t\geq 0$.
\end{theorem}

As we will see in Theorem \ref{qqq}, Ricci flows with unbounded curvature at later times were first constructed in higher dimensions in the context of flows with nonnegative complex sectional curvature.

To see a second way in which the classical Hamilton-Shi flow cannot be a substitute for the flow from Theorem \ref{fullthm}, consider for a moment the restricted situation that $(\mm,g_0)$ is both complete and of bounded curvature. We can flow such a metric not only with Theorem \ref{fullthm} but also with the Hamilton-Shi flow. By the uniqueness of complete Ricci flows asserted in Theorem \ref{fullthm}, they must agree -- at least while the Hamilton-Shi flow exists. By Theorem \ref{hamilton_shi}, the Hamilton-Shi flow exists for all time unless the curvature blows up, so one might naively think that this should be when the flow from Theorem \ref{fullthm} stops too. However, 
this more general flow can typically keep going.

\begin{theorem}[{\bf Flowing beyond curvature blow-up;} \cite{GT4}]
There exists a complete surface $(\mm,g_0)$ of bounded curvature such that the subsequent Ricci flow $g(t)$ given by Theorem \ref{fullthm} exists for all time $t\geq 0$, and that for some times $0<\tilde T<t_0<\infty$ we have:
\begin{itemize}
\item
The curvature of $g(t)$ is unbounded as $t\upto \tilde T$, but bounded on $[0,T_0]$ for all $T_0\in [0,\tilde T)$.
\item
The curvature is unbounded for each $t\geq t_0$.
\end{itemize}
\end{theorem}

Clearly, the Hamilton-Shi flow will agree with $g(t)$ above until time $\tilde T$, when it stops, leaving $g(t)$ to continue alone.

\begin{figure}
\begin{tikzpicture}[xscale=2.3,yscale=2.3]
\draw [ultra thick, <->] (0,1.3) node [left] {$\sup_\mm |K|$} -- (0,0) -- (4.2,0) node [below] {$t$};

%\draw[ultra thick, domain=0:0.8] 
\draw[line width = 2.0pt, domain=0:0.8] 
plot (\x, {
0.3/(1-\x)
});

%\draw[ultra thick, domain=2.2:4] 
\draw[line width = 2.0pt, domain=2.2:4] 
plot (\x, {
0.3/(\x-2)
});

% cheating with the 0.9 - should be 1
\draw [thick] (0.9,-.02) node[below, align=center]
{$\tilde T$
} -- (0.9,0.02);
\draw [thick, dashed] (0.9,0) -- (0.9,1.5);

\draw [thick] (2.1,-.02) node[below]{$t_2$} -- (2.1,0.02);
\draw [thick, dashed] (2.1,0) -- (2.1,1.5);

\draw[thick, <->] (1.1,1) -- (1.5,1) node[below, align=center]{curvature\\ unbounded} --(1.9,1);

\draw [thick] (1.1,-.02) node[below]{$t_0$} -- (1.1,0.02);
\draw [thick] (1.9,-.02) node[below]{$t_1$} -- (1.9,0.02);

\end{tikzpicture}
\caption{The flow in Theorem \ref{bursts_thm} switches back and forth between bounded and unbounded curvature, remaining always complete. The Hamilton-Shi flow ends at $t=\tilde T$.}
\label{bursts_fig}
\end{figure}
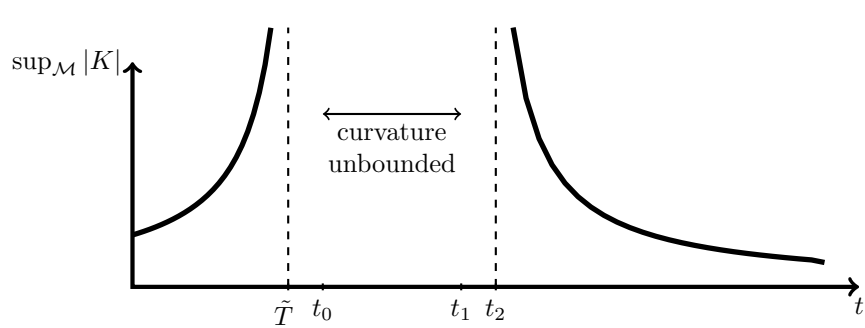

An earlier result of Cabezas-Rivas and Wilking that we state in Theorem \ref{qqq} constructs examples of Ricci flows in dimensions four and higher, with nonnegative complex sectional curvature, with unbounded curvature precisely for $t\geq 1$.
In fact, more exotic behaviour can be engineered in which the flow passes back and forth between periods in which the curvature of the flow is unbounded and bounded, such as in the following theorem and Figure \ref{bursts_fig}.

\begin{theorem}[{\bf Bursts of unbounded curvature;} \cite{GT4}]
\label{bursts_thm}
There exists a complete surface $(\mm,g_0)$ of bounded curvature such that the subsequent Ricci flow $g(t)$ given by Theorem \ref{fullthm} exists for all time $t\geq 0$, and such that for some times $0<\tilde T<t_0<t_1<t_2<\infty$ we have:
\begin{itemize}
\item
The curvature of $g(t)$ is bounded on $[0,T_0]$ for all $T_0\in [0,\tilde T)$, but unbounded on $[0,\tilde T)$.
\item
The curvature is unbounded for each $t\in [t_0,t_1]$.
\item
The curvature is bounded for $t\geq t_2$.
\end{itemize}
\end{theorem}

%\vskip 0.1cm

Of course, by uniqueness, the only lever we have to engineer this behaviour is the choice of the initial metric. Once the flow starts, it is determined forever.

Now we have a Ricci flow starting with an arbitrary initial surface, we can return to the question of whether this Ricci flow will geometrise the surface, finding a metric of constant curvature. If we restrict our discussion to initial surfaces $(\mm,g_0)$ that are conformal to a hyperbolic metric $H$, then by definition of $T$ in Theorem \ref{fullthm}, our Ricci flow must exist for all $t\geq 0$ and we can ask the question as to whether the flow will manage to converge to $H$.
We find that it does:
\begin{theorem}[{\bf Ricci flow geometrises general hyperbolic surfaces;} \cite{GT2}]
\label{asymp_thm}
If $(\mm,g_0)$ has a conformally equivalent, complete, hyperbolic metric $H$, then the Ricci flow $g(t)$ from Theorem \ref{fullthm} exists for all $t\geq 0$ and finds $H$ in the sense that
$$\frac{g(t)}{2t}\to H$$
smoothly locally as $t\to\infty$.
\end{theorem}

We view this type of result not as a route to giving an alternative proof of the Uniformisation theorem for general surfaces, but partly as a stepping stone to proving higher-dimensional uniformisation results, 
partly as a way of comparing the geometry of general metrics to that of conformal constant curvature metrics -- for example the determinant of the Laplacian behaves well under Ricci flow 
\cite{AAR} -- 
and partly as a means for understanding how Ricci flow organises itself in order to find a special metric.
In the latter direction, it is interesting to apply Theorem \ref{asymp_thm} in the case that $\mm$ is noncompact (for example, the open disc) and $g_0$ is a metric supplied from Theorem \ref{unbdd_ex}. In this case we know on one hand that $\frac{g(t)}{2t}$ converges to a hyperbolic metric, i.e. one of constant curvature $-1$, and on the other hand that it remains with unbounded curvature for all time.
These two assertions are not contradictory because the convergence to a hyperbolic metric is smooth \emph{local} convergence. In other words, the Ricci flow cannot prevent unbounded curvature, but it does organise itself to force regions of large curvature out towards spatial infinity.

\section{Flows with unbounded curvature in higher dimensions}

\begin{figure}[h]
\begin{center}
\def\svgwidth{300pt}
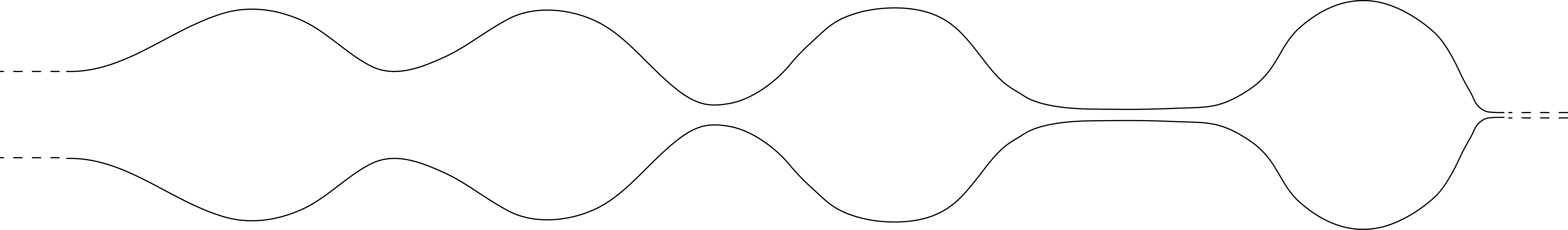
\end{center}
\caption{$S^2\times \RR$ with thinner and thinner necks}
%%ALWAYS PUT LABEL AFTER CAPTION!!!!
\label{bubble_fig}
\end{figure}

It is interesting to speculate as to how many of the results of the last section generalise to higher dimensions. Certainly the existence theory in Theorem \ref{fullthm} cannot possibly be expected to hold in such generality because one can choose smooth, complete, three-dimensional Riemannian manifolds of unbounded curvature that we cannot hope to evolve under Ricci flow, even for a short time.
For example, one could take the underlying manifold to be $S^2\times\RR$ and endow it with a warped product metric so that metrically it consists of an infinite chain of three-spheres connected by thinner and thinner (and longer and longer) necks, as in Figure \ref{bubble_fig}.
However small we take $\ep>0$, if we pick a neck that is sufficiently thin and long, then the Ricci flow will be inclined to pinch it 
within time $\ep$. (See \cite[\S 1.3.2]{RFnotes} for more about this type of neck-pinch singularity.)

One can, however, hope to flow manifolds of unbounded curvature that also satisfy an appropriate nonnegative curvature condition that can rule out neck-pinch type singularities happening in an arbitrarily short time. An interesting scenario would be to consider complete three-dimensional manifolds of nonnegative Ricci curvature. It is very likely that Ricci flow knows how to flow such manifolds, preserving the nonnegativity of the Ricci curvature, in a unique way. 

In 1989, Shi  \cite{shi_nonnegRic} classified complete three manifolds of nonnegative Ricci curvature that have \emph{bounded curvature}. The key step, following Hamilton \cite{ham4PCO}, was to flow forwards using the Ricci flow and prove that either the Ricci curvature becomes strictly positive immediately, or the manifold splits locally as a product. Being able to drop the bounded curvature assumption while still being able to flow even for a short time would have 
%immediately 
yielded a classification of \emph{all} three manifolds of nonnegative Ricci curvature, which is the most natural class to consider. This classification had to wait over twenty years, and the development of alternative minimal surface techniques by Liu \cite{gangliu}.
The natural Ricci flow question of existence with unbounded curvature remains open.

A very natural situation in which the existence side of the theory \emph{has} been successfully developed is the case that we flow manifolds of nonnegative \emph{complex} sectional curvature. 
This curvature condition implies nonnegative sectional curvature, and is implied by nonnegative curvature operator; it is preserved under Ricci flow \cite{CR_W}.
Cabezas-Rivas and Wilking proved existence of solutions in this situation, retaining the curvature condition.
\begin{theorem}[Cabezas-Rivas and Wilking {\cite[Theorem 1]{CR_W}}]
\label{CRW1}
Given any complete Riemannian manifold $(\mm,g_0)$ with nonnegative complex sectional curvature, there exists a Ricci flow $g(t)$ for $t\in [0,\ep)$, some $\ep>0$, with $g(0)=g_0$ and with $g(t)$ having nonnegative complex sectional curvature for each $t\in [0,\ep)$.
\end{theorem}
It is an 
interesting open question to prove uniqueness of this solution -- conceivably there could be many other solutions with the same initial metric. 
Indeed, \emph{a priori} the Ricci flow could have infinitely many nonunique branches at each instant of time, as is the case for parabolic equations on bounded domains.

Except in low dimensions, the restriction of nonnegative complex sectional curvature does not in itself enforce boundedness of the curvature of a Ricci flow as time advances:

\begin{theorem}[Cabezas-Rivas and Wilking {\cite[Theorem 4, Corollary 3]{CR_W}}]
\label{qqq}
There exist Ricci flows with nonnegative complex sectional curvature for $t\in [0,\infty)$ with unbounded curvature for all time, and others with unbounded curvature precisely for $t\geq 1$.
On the other hand, if for an $n$-dimensional Riemannian manifold $(\mm,g_0)$ with nonnegative complex sectional curvature there exists $v_0>0$ such that
$$\Vol_{g_0}(B_{g_0}(x,1))\geq v_0$$
for all $x\in\mm$, then we may assume that the curvature of the Ricci flow arising in Theorem \ref{CRW1} is bounded above by $\frac{C(n,v_0)}{t}$.
\end{theorem}

The existence theory above is not explicit about the length of time for which the solution will persist. However, by virtue of the nonnegativity of the complex sectional curvature, a lower bound for the existence time can be read off from the initial geometry of the flow:

\begin{theorem}[Cabezas-Rivas and Wilking {\cite[special case of Corollary 5]{CR_W}}]
Given an $n$-dimensional manifold $\mm$, 
there exists $\ep>0$ depending only on $n$ such that
we may assume that a Ricci flow arising from Theorem \ref{CRW1}
exists for $t\in [0,T)$ where
$$T=\ep.\sup\left\{\frac{\Vol_{g_0}(B_{g_0}(x,r))}{r^{n-2}}\ \bigg|\ x\in\mm, r>0\right\}\in (0,\infty].$$
\end{theorem}

Finally we remark that perhaps the most natural context in which one might hope to generalise the results of Section \ref{2Dsect} is that of higher dimensional \emph{K\"ahler} Ricci flow. We leave this discussion for another occasion.

\section{Flowing rough metrics or Alexandrov spaces}

The entire discussion so far has considered only Ricci flows starting with \emph{smooth} Riemannian metrics. One can also ask whether it is possible to start flowing from a metric that is not smooth, or is possibly not even a Riemannian metric at all, but perhaps a metric space with some basic structure. A Ricci flow arising in this way would generally have unbounded curvature at least in the limit $t\downto 0$.

One might naively think that if we can view Ricci flow as a parabolic equation, then it should be irrelevant what the regularity of an initial Riemannian metric is because Ricci flow should instantly smooth it out. However, Ricci flow will \emph{not} smooth in the same quantified way as the ordinary heat equation. The Harnack inequality, in the sense of \cite[Section 1]{harnack_survey}, can be used to prove that positive solutions of the heat equation on, say, Euclidean space must smear out at a certain rate, with the heat kernel being the extreme case. 
By contrast, the coefficient $e^{-2u}$ in \eqref{2DRFeq} can interrupt this behaviour, as is exploited implicitly in the proofs of Theorems \ref{unbdd_ex} and \ref{bursts_thm}. This effect from a PDE perspective says that a delta-mass for $e^{2u}$ will remain a delta-mass for a while under the evolution equation \eqref{2DRFeq}, and will not spread out like a heat kernel. See the discussion in \cite[Section 8.2]{vaz}.

Nevertheless, there are several situations in which it is possible to start the Ricci flow with a very rough object, perhaps a certain type of metric space.
The general principle behind most results of this form, as well as some of the results we have discussed earlier in this survey, is that one approximates the initial data by a sequence of smooth Riemannian manifolds $(\mm_i,g_i)$, flows each of these, and then tries to take a limit of these smooth Ricci flows. The real art is to prove uniform estimates on the sequence of smooth Ricci flows and their existence time so that this limit can be taken. A number of estimates of this form could be summarised loosely as
\begin{quote}
\fbox{
\parbox{140\unitlength}{
Initial metric $g_0$ has\\ 
(i) lower curvature bound, and
(ii) noncollapsing hypothesis
}
}
$\implies$
\fbox{
\parbox{120\unitlength}{
Curvature decays like\\ 
\makebox[120\unitlength][c]{
$|\Rm_{g(t)}|\leq \frac{C}{t}$
}
}
}
\end{quote}
Without the noncollapsing condition (ii), no such uniform estimate can hold as we see by returning to the `contracting cusp' example of Theorem \ref{revcusp_special} where the curvature decay is expected to be like $C/t^2$.

Although we do not attempt a complete survey of such results, we do wish to highlight one of the results of M. Simon \cite{simon2012} of this form.
\begin{theorem}[Special case of {\cite[Theorem 7.1]{simon2012}}]
\label{simon}
Given a closed three or two-dimensional manifold 
$(\mm,g_0)$ that satisfies 
\begin{enumerate}[(a)]
\item
$\Ric_{g_0}\geq k\in \RR$, and 
\item
$\Vol_{g_0}(B_{g_0}(x,1))\geq v_0>0$ for all $x\in\mm$,
\end{enumerate}
there exist constants $T=T(k,v_0)>0$ and $K=K(k,v_0)>0$ such that the Ricci flow $g(t)$ with $g(0)=g_0$ exists for $t\in [0,T)$ and 
satisfies, for each $t\in (0,T)$, the inequalities
\begin{enumerate}[(i)]
\item
$\Ric_{g(t)}\geq -K$, 
\item
$\Vol_{g(t)}(B_{g(t)}(x,1))\geq \frac{v_0}{2}$ for all $x\in\mm$,
\item
$|\Rm_{g(t)}|\leq \frac{K}{t}$, and
\item
$e^{K(t-s)}d_{g(s)}(p,q)\geq d_{g(t)}(p,q)\geq d_{g(s)}(p,q)-K(\sqrt{t}-\sqrt{s})$ for all $s\in (0,t]$ and $p,q\in\mm$, where
$d_{g(s)}(p,q)$ is the Riemannian distance between $p$ and $q$ with respect to $g(s)$.
%$e^{K(t-s)}d(p,q,s)\geq d(p,q,t)\geq d(p,q,s)-K(\sqrt{t}-\sqrt{s})$ for all $s\in (0,t]$ and $p,q\in\mm$.
\end{enumerate}
\end{theorem}

As a consequence, Simon \cite{simon2012} was able to run the 
Ricci flow for a definite amount of time, starting with any metric space arising as a Gromov-Hausdorff limit of smooth closed three-manifolds satisfying (a) and (b) of Theorem \ref{simon}
uniformly.

Of course, we would like to be able to describe \emph{synthetically} the metric spaces that can be flowed, rather than describe them as limit points under Gromov-Hausdorff convergence.
A clean situation  in which this can be done was described by 
T. Richard \cite{TR}. The starting point for his work is the notion of Alexandrov space, which for our purposes is a type of metric space $(X,d)$ that satisfies a very weak notion of curvature bounded below. Such spaces automatically have integral (or infinite) Hausdorff dimension, and we require this dimension to be two. The resulting object turns out to be  topologically a surface (possibly with boundary), and we require this surface to be closed.
These constraints define the notion of \emph{Alexandrov surface} -- see 
\cite{richard_thesis} and \cite{BBI} for a more precise definition and further details.

Following Alexandrov, Richard \cite{richard_thesis, TR} proved that Alexandrov surfaces can be approximated by smooth Riemannian surfaces satisfying (a) and (b) of Theorem \ref{simon}
uniformly, and hence showed that Ricci flow smooths out an Alexandrov surface. He also proved that it does this in a unique way.
\begin{theorem}[{\bf Flowing Alexandrov surfaces} \cite{TR}]
Given an Alexandrov surface $(X,d)$ with curvature bounded below by 
$-1$, there exists a Ricci flow $g(t)$ on a closed surface $\mm$, for $t\in (0,T)$, some $T>0$, with $K_{g(t)}\geq -1$ and such that the distance function $d_{g(t)}:\mm\times\mm\to [0,\infty)$ of $g(t)$ converges uniformly to some distance function $d_0:\mm\times\mm\to [0,\infty)$ as $t\downto 0$, where $(\mm,d_0)$ is isometric to $(X,d)$.

Moreover, the Ricci flow is unique up to isometry in the sense that if $(\hat\mm,\hat g(t))$, defined for $t\in (0,\hat T)$,
is any other Ricci flow with curvature uniformly bounded from below, 
for which
$d_{\hat g(t)}\to \hat d_0$ as $t\downto 0$ with 
$(\hat\mm,\hat d_0)$ isometric to $(X,d)$, then there exists a smooth map $\vph:\mm\to\hat\mm$ that is an isometry from
$(\mm,g(t))$ to $(\hat\mm,\hat g(t))$ for each $t\in(0,\min\{T,\hat T\})$.
\end{theorem}
A curious consequence of this theorem is that Ricci flow knows again how to organise itself, this time to uniquely and instantaneously endow such metric spaces 
with a conformal structure.

%\frenchspacing

\end{document}